\documentstyle[12pt]{article}

\def\A{{\cal A}}

\def\H{{\cal H}}
\def\h{{\bf h}}
\def\k{{\kappa}}
\def\n{{\nu}}
\def\L{{\cal L}}
\def\M{{\cal M}}
\def\RR{{\bf R}}
\def\bw{{\overline{w}}}
\def\w{{\rm w}}
\def\bz{{\overline{z}}}
\def\z{{\zeta}}
\def\sa{{\rm sa}}
\def\Re{{\rm Re}\>}
\def\Im{{\rm Im}\>}
\def\ip<#1|#2>{\left<#1\vphantom{#2}\right|\left.\vphantom{#1}#2\right>}
\def\<#1|#2>{\left<\right.#1\left|\right.#2\left.\right>}
\def\(#1){\left(#1\right)}
\def\[#1|#2]{\left\lbrace #1\vphantom{#2}\right|
             \left.\vphantom{#1}#2\right\rbrace}

\begin{document}
\date{March 1995}
\title{Endomorphism Semigroups and Lightlike Translations}
\author{D.~R.~Davidson*\\
Dipartimento di Matematica\\
Universit\`a di Roma ``La Sapienza''$\quad$00185 Roma, Italy\\
Email:  davidson@mat.uniroma1.it
}
\maketitle
\begin{abstract}
Borchers and Wiesbrock have studied the one-parameter semigroups of
endomorphisms of von Neumann algebras that appear as lightlike translations in
the theory of algebras of local observables, showing that they automatically
transform under the appropriate modular automorphisms as under velocity
transformations.  These results are here abstracted and analyzed as essentially
operator-theoretic.  Criteria are then established for a spatial derivation of
a von Neumann algebra to generate a one-parameter semigroup of endomorphisms,
and all of this is combined to establish a von Neumann-algebraic converse to
the Borchers and Wiesbrock results.  This sort of analysis is then applied to
questions of isotony and covariance for local algebras, to show that Poincar\'e
covariance together with a domain condition for the translations can imply
isotony.

\bigskip\noindent
* This research was partly supported by a fellowship from the
Consiglio Nazionale delle Ricerche.
\end{abstract}
\newpage

\section{Introduction}

The standard situation for a pair of complementary spacetime regions in the
theory of algebras of local observables, under the assumption of duality in the
vacuum sector, is just that termed standard in the theory of von Neumann
algebras:  we have a von Neumann algebra $\M$ and its commutant $\M'$ acting on
a Hilbert space $\H$, with a common cyclic and separating unit vector $\Omega$,
the vacuum vector.  In the particular situation in which $\M$ and $\M'$
correspond to the algebras of observables for a pair of complementary wedge
regions (for definiteness let us take them to be $W_R=\[x|x_1>{|t|}]$ and
$W_L=\[x|x_1<-{|t|}]$ respectively) it is expected that the modular
automorphism group $\sigma_t(A)=\Delta^{it}A\Delta^{-it}$ will correspond to
the Lorentz velocity transformations $V_1(2\pi t)$ in the direction $\hat x_1$
orthogonal to the vertex $x_1=t=0$ of the wedges, and that the modular
conjugation $J$ will be a slight variant of the TCP operator $\Theta$, namely
$J=\Theta R$ where $R$ is a rotation by the angle $\pi$ about the direction
$\hat x_1$ \cite{BW}.  In that case the lightlike translations
$U(a)=T(a(\hat x_0+\hat x_1))$ will be a strongly continuous one-parameter
group of unitary operators on $\H$, which should have the following four
properties:\\
(a) By Lorentz and TCP covariance,
$\Delta^{it}U(a)\Delta^{-it}=U(e^{-2\pi t}a)$ and $JU(a)J=U(-a)$ for all
real $a$ and $t$;\\
(b) By the spectral condition,
$U(a)$ should have a positive generator $H$;\\
(c) By isotony,
for $a\geq 0$ the corresponding adjoint action $A\rightarrow U(a)AU(-a)$ should
give a one-parameter semigroup of endomorphisms of $\M$ (and thus for $a\leq 0$
likewise of $\M'$); and, finally,\\
(d) The vacuum vector $\Omega$ should be fixed by all $U(a)$, and thus
annihilated by $H$.

In this connection Borchers has shown that these four conditions are not all
independent:  in particular, if the last three hold, then the covariance
conditions follow automatically \cite{Bo}.  Wiesbrock then proved conversely
that if (a), (c), and (d) hold, $U(a)$ automatically has a positive generator
\cite{Wb}.  In this note we first analyze the results of Borchers and of
Wiesbrock, showing that they are essentially operator-theoretic statements
about relations between $J$, $\Delta$, and $H$; we then demonstrate that they
are part of a larger chain of converses, in which we separate out the von
Neumann-algebraic content, and in the process perhaps shed some further light
on these remarkable theorems.  Specifically, we show that if the generator $H$
gives a derivation $\delta$ of $\M$ satisfying certain additional conditions,
then any three of the conditions listed above for $U(a)$ together imply the
fourth.  Note that in the local algebra context, it can be shown that $\M$ and
$\M'$ must be Type III${}_1$ factors \cite{Dr}, but this will not be used in
the following; the results will simply be stated in terms of arbitrary von
Neumann algebras.  We then proceed to apply these methods to questions of
isotony and covariance in the context of the theory of local algebras.

The situation here is analogous to, but in some respects altogether different
from, the case of spatial derivations that generate automorphism groups of
von Neumann algebras, which has been extensively studied
(\cite{OAQSM}, Section 3.2.5, and references therein).  We will develop the
analogy more specifically in Section III, but the relevant condition in the
automorphism case is that $U(a)$ should commute with $J$ and with all
$\Delta^{it}$; then the key question is to determine precisely what additional
conditions on a derivation $\delta$ suffice to show that it generates an
automorphism group.  The best result in this direction is that of \cite{BH}, in
which the only additional assumption is that the derivation has a domain
$D(\delta)$ such that $D(\delta)\Omega$ is a core for $H$.  However, the proof
of this result is rather difficult, and does not generalize to the endomorphism
case.  We will generally make do with more restrictive conditions here, but it
would be interesting to determine precisely what conditions suffice to
guarantee that $\delta$ generates an endomorphism semigroup; in Section V
we discuss a case in which a condition like that of \cite{BH} suffices.

In what follows, Section II presents the results of Borchers and Wiesbrock,
simplified and reduced to their operator-theoretic essence; Section III
discusses criteria for a spatial derivation to generate a one-parameter
semigroup of endomorphisms; Section IV collects the resulting information
about lightlike translations; and Section V discusses applications of this
analysis to certain systems of local algebras, showing that isotony relations
hold given only covariance and a core condition as above for the translations.

\section{Commutation Relations}

The following is of course based heavily on \cite{Bo}, with contributions from
\cite{Wb} and \cite{BGL}.  We present these results in what would seem
to be their natural setting:  that of commutation relations for one-parameter
unitary groups, and their unbounded generators.

\bigskip
\noindent{\bf Theorem 1:}~~
{\em
Let $V(\lambda)=\Delta^{i\lambda/2\pi}$ and $U(a)=e^{iaH}$ be two strongly
continuous one-parameter unitary groups.  Then any two of the following
conditions imply the third:\\
(a) $V(\lambda)U(a)V(-\lambda)=U(e^{-\lambda}a)$ for all real $a$ and
$\lambda$;\\
(b) $H$ is positive;\\
(c) $\Delta^{1/2}U(a)\supset U(-a)\Delta^{1/2}$ for all $a\geq 0$.
}

\smallskip\goodbreak
\noindent{\bf Proof of Theorem 1:}~~
{
Here $\Delta$ and $H$ are unbounded operators with their natural domains of
definition, $\Delta$ positive and $H$ self-adjoint by Stone's Theorem; there
is therefore a dense set $D_\omega$ of vectors $\psi$ for which
$\Delta^{iz}\psi=V(2\pi z)\psi$ is entire analytic.  Let us take two fixed
vectors $\psi,\phi\in D_\omega$, and define two functions
\begin{equation}
F(z,w)=\ip<\Delta^{i\bz}\psi|U(e^{2\pi w})\Delta^{iz}\phi>\quad\hbox{\rm and }
\quad G(z,w)=\ip<\Delta^{i\bz}\psi|U(-e^{2\pi w})\Delta^{iz}\phi>,
\end{equation}
both entire analytic functions of $z$ for real $w$.  If $H$ is positive, then
$F$ and $G$ will have jointly analytic continuations, continuous at the
boundary, to $0\leq\Im w\leq 1/2$ and $-1/2\leq\Im w\leq 0$ respectively,
satisfying
\begin{equation}
|F(z,w)|\leq\|\Delta^{\Im z}\psi\|\>\|\Delta^{-\Im z}\phi\|\quad\hbox{\rm and }
\quad|G(z,w)|\leq\|\Delta^{\Im z}\psi\|\>\|\Delta^{-\Im z}\phi\|
\end{equation}
independent of $w$ and $\Re z$ over the whole region of definition; these
bounds hold for $\Im w=0$ whether $H$ is positive or not.

If (a) and (b) hold, then $F$ satisfies the complex identity
$F(z,w)=F(0,z+w)$.  Taking $w$ real, $z=i/2$, this implies that
$\<\psi|U(-a)\phi>=\<\Delta^{1/2}\psi|U(a)\Delta^{-1/2}\phi>$ for all
$a\geq 0$.  Since $\phi'=\Delta^{-1/2}\phi\in D_\omega$ if and only if
$\phi\in D_\omega$, we may equally write
\begin{equation}
\ip<\psi|U(-a)\Delta^{1/2}\phi'>=\ip<\Delta^{1/2}\psi|U(a)\phi'>.
\end{equation}
Since $\psi$ may be any vector in a dense set, this in fact holds for any
$\psi\in D(\Delta^{1/2})$.  The right-hand side is a bounded function of
$\phi'$, so that the left-hand side is also; this implies that
$U(a)\psi\in D(\Delta^{1/2})$, and we may write
\begin{equation}
\ip<\Delta^{1/2}U(a)\psi|\phi'>=\ip<U(-a)\Delta^{1/2}\psi|\phi'>.
\end{equation}
Since $\phi'$ may be any vector in a dense set,
$\Delta^{1/2}U(a)\psi=U(-a)\Delta^{1/2}\psi$ for every
$\psi\in D(\Delta^{1/2})$ and any $a\geq 0$; but this is just the
statement of (c).

On the other hand, if we are given (c), then we have for real $t,s$
\begin{eqnarray}
F(t+i/2,s)&=&\ip<\Delta^{1/2}\Delta^{it}\psi
|U(e^{2\pi s})\Delta^{-1/2}\Delta^{it}\phi>\nonumber\\
&=&\ip<\psi|\Delta^{-it}U(-e^{2\pi s})\Delta^{1/2}\Delta^{-1/2}\Delta^{it}\phi>
=G(t,s).
\end{eqnarray}
Not only is $|F(t,s)|\leq\|\psi\|\|\phi\|$ for real $t,s$, but also
$|F(t+i/2,s)|\leq\|\psi\|\|\phi\|$.   In addition, the bound given above for
$\Im w=0$, depending only on $\Im z$, implies that $F(z,w)$ is bounded on the
strip $0\leq\Im z\leq 1/2$, $\Im w=0$; then by Hadamard's three-line theorem,
$|F(z,w)|\leq\|\psi\|\|\phi\|$ on this strip.  If in addition (a) holds, then
the identity $F(z,w)=F(z+w,0)$ shows that $|F(z,w)|\leq\|\psi\|\|\phi\|$ for
$0\leq\Im z+\Im w\leq 1/2$.  Thus
$|\<\psi|e^{-H}\phi>|=|F(0,i/4)|\leq\|\psi\|\|\phi\|$; but $\psi,\phi$ may be
any vectors in a dense set, so $\|e^{-H}\|\leq 1$ and $H$ is positive.

If instead (b) and (c) hold, then also $F(t,s+i/2)=G(t,s)$ for real $t,s$.
 From
\begin{equation}
F(t,s+i/2)=G(t,s)=F(t+i/2,s)
\end{equation}
it follows that $F(z,w)$ and $G(z-i/2,w)$ have a common jointly entire analytic
continuation $\tilde F(z,w)$ satisfying $\tilde F(z,w+i/2)=\tilde F(z+i/2,w)$
for all $z,w$.  Consider $\tilde F(t+\z,s-\z)$ as a function of $\z$:  it is
periodic with period $i/2$, and its modulus has an upper bound depending only
on $\Im\z$; therefore it is bounded, hence constant.  It follows that
$\tilde F(z,w+\z)=\tilde F(z+\z,w)$ for all complex $\z$.  In particular,
$F(t,s)=F(0,t+s)$, so that
$\<\psi|V(\lambda)U(a)V(-\lambda)\phi>=\<\psi|U(e^{-\lambda}a)\phi>$ for all
real $\lambda$ and all $a\geq 0$.  Since $\psi,\phi$ may be any vectors in a
dense set, $V(\lambda)U(a)V(-\lambda)=U(e^{-\lambda}a)$ for all real
$\lambda$ and all $a\geq 0$, and taking adjoints we obtain the result for
all real $a$.
}
\smallskip

Examining the proof, we see that (c) is in fact equivalent to a number of other
statements, a suitably weak one being for example that for all $a\geq 0$,
$\Delta^{1/2}U(a)$ is equal to $U(-a)\Delta^{1/2}$ in the sense of quadratic
forms, with form domain $D_\omega$.
\goodbreak

The conditions of Theorem 1 are trivially satisfied if $H=0$ no matter what
$V(\lambda)$ may be; furthermore, if $\psi$ is a vector invariant under all
$U(a)$, then so also is $V(\lambda)\psi$ for any $\lambda$.  The Hilbert space
can thus be decomposed into a direct sum of the subspace of such vectors, and
its orthogonal complement, on which $V(\lambda)$ and $U(a)$ are non-trivial.
Then the following is a simple computation:

\bigskip\goodbreak
\noindent{\bf Proposition 2:}~~
{\em
Given a representation of the relations of Theorem 1, if $H$ is restricted
to the orthogonal complement of its null space, then it has a self-adjoint
logarithm $P$, such that $H$ so restricted equals $T(i)$ where $T(x)=e^{-ixP}$.
This satisfies
\begin{equation}
V(\lambda)T(x)=e^{i\lambda x}T(x)V(\lambda)
\end{equation}
for all real $x$ and $\lambda$, so that $V(\lambda)$ and $T(x)$ give a
representation of the canonical commutation relations in Weyl form.
Conversely,
any representation of these commutation relations gives a
representation of the relations of Theorem 1 in this manner.
}
\bigskip\goodbreak

Thus the Stone-von Neumann classification of representations of the canonical
commutation relations provides a classification of representations of the
relations in Theorem 1:  up to multiplicity, the only representations are
either trivial (with $V(\lambda)$ arbitrary), or else given by
$\Delta=e^{2\pi X}$ and $H=e^P$, where $X$ and $P$ have the familiar
Schr\"odinger form for the commutation relations $[X,P]=i$.

To adapt these results to modular theory, we need to extend them slightly
so as to include the modular conjugation.

\bigskip\goodbreak
\noindent{\bf Theorem 3:}~~
{\em
In addition to the premises of Theorem 1, let $J$ be a complex conjugation
(an antiunitary involution) commuting with all $V(\lambda)$.  Then any two of
the following conditions imply the third:\\
(a${}'$) $V(\lambda)U(a)V(-\lambda)=U(e^{-\lambda}a)$ and $JU(a)J=U(-a)$ for
all real $a$ and $\lambda$;\\
(b) $H$ is positive;\\
(c${}'$) $J\Delta^{1/2}U(a)\supset U(a)J\Delta^{1/2}$ for all $a\geq 0$.\\
To show (b) it is not necessary to assume that $JU(a)J=U(-a)$.
}

\smallskip\goodbreak
\noindent{\bf Proof of Theorem 3:}~~
{
The proof is quite similar to that of Theorem 1.  If (a${}'$) and (b) hold,
then by Theorem 1, we have
$J\Delta^{1/2}U(a)\supset JU(-a)\Delta^{1/2}=U(a)J\Delta^{1/2}$ for all
$a\geq 0$, and we are finished immediately.  For the remaining parts we
define two additional functions with properties like those of $G$ and $F$,
\begin{equation}
H(z,w)=\ip<\Delta^{i\bz}\psi|JU(e^{2\pi\bw})J\Delta^{iz}\phi>
\>\,\hbox{\rm and }\>
K(z,w)=\ip<\Delta^{i\bz}\psi|JU(-e^{2\pi\bw})J\Delta^{iz}\phi>.
\end{equation}

If (c${}'$) holds, then we also have the adjoint statement
$\Delta^{1/2}JU(-a)\supset U(-a)\Delta^{1/2}J$ for all $a\geq 0$.
Instead of one identity, we now have two:
\begin{eqnarray}
F(t+i/2,s)&=&\ip<\Delta^{1/2}J^2\Delta^{it}\psi|U(e^{2\pi s})
\Delta^{-1/2}\Delta^{it}\phi>\nonumber\\
&=&\ip<\psi|\Delta^{-it}JU(e^{2\pi s})J\Delta^{1/2}
\Delta^{-1/2}\Delta^{it}\phi>=H(t,s);\\
K(t+i/2,s)&=&\ip<\Delta^{1/2}\Delta^{it}\psi|JU(-e^{2\pi s})J
\Delta^{-1/2}\Delta^{it}\phi>\nonumber\\
&=&\ip<\psi|\Delta^{-it}U(-e^{2\pi s})\Delta^{1/2}
\Delta^{-1/2}\Delta^{it}\phi>=G(t,s).
\end{eqnarray}
If in addition (a) of Theorem 1 holds, then the argument of Theorem 1 applies
(with the substitution of $H(z,w)$ for $G(z,w)$) to show that $H$ is positive.

If instead (b) and (c${}'$) hold, then the analytic continuation argument is
only slightly more complicated.  We have already $F(t+i/2,s)=H(t,s)$
and $K(t+i/2,s)=G(t,s)$; in addition now
$F(t,s+i/2)=G(t,s)$ and $K(t,s+i/2)=H(t,s)$.
It follows that $F(z,w)$, $G(z,w-i/2)$, $H(z-i/2,w)$ and $K(z-i/2,w+i/2)$ have
a common jointly entire analytic continuation $\tilde F(z,w)$ satisfying
$\tilde F(z,w+i)=\tilde F(z+i,w)$ for all $z,w$.

For real $t,s$, the function $\tilde F(t+\z,s-\z)$ of $\z$ again has a bound
independent of $\Re\z$, and now it is periodic with period $i$; therefore it is
again bounded and constant.  Thus $\tilde F(z,w+\z)=\tilde F(z+\z,w)$ for all
complex $\z$, and as in Theorem 1, $V(\lambda)U(a)V(-\lambda)=U(e^{-\lambda})$
for all real $\lambda$ and $a$.  Furthermore $K(t,s)=F(t,s)$, so that
$\<\psi|JU(-a)J\phi>=\<\psi|U(a)\phi>$ for all $a\geq 0$.  Since $\psi,\phi$
may be any vectors in a dense set, $JU(-a)J=U(a)$ for all $a\geq 0$, and taking
adjoints we obtain the result for all real $a$.
}
\smallskip

Proposition 2 still holds for the representations of the relations of
Theorem 3, but now $J$ is a complex conjugation such that
$JV(\lambda)J=V(\lambda)$ and $JT(x)J=T(-x)$.  Thus in an irreducible
Schr\"odinger representation $J$ is the complex conjugation in $P$-space
(up to multiplication or conjugation by a complex phase).
\goodbreak

\section{Endomorphism Semigroups}

If we have a von Neumann algebra $\M$ and its commutant $\M'$ acting on a
Hilbert space $\H$, with a common cyclic and separating vector $\Omega$,
we may define real linear spaces $R=\overline{\M^\sa\Omega}$ and
$R'=\overline{\M^{'\sa}\Omega}$.  Then $\<\psi|\phi>$ is real for all
$\psi\in R,\phi\in R'$, and furthermore $R'$ is precisely the set of all
$\psi$ such that $\<\psi|\phi>$ is real for all $\phi\in R$.  Also, both
$D(\Delta^{1/2})=R+iR$ and $D(\Delta^{-1/2})=R'+iR'$ are dense in $\H$,
and we have $R=\[\psi|\psi\in D(\Delta^{1/2}), J\Delta^{1/2}\psi=\psi]$ and
$R'=\[\psi|\psi\in D(\Delta^{-1/2}), J\Delta^{-1/2}\psi=\psi]$.
Thus the condition $J\Delta^{1/2}U(a)\supset U(a)J\Delta^{1/2}$ corresponds
to $U(a)R\subset R$.

For any $\psi\in R$, there is a sequence $X_n\in\M^\sa$ such that
$X_n\Omega\rightarrow\psi$, but there need not be a bounded operator
$X\in\M^\sa$ such that $X\Omega=\psi$; in general there is only a closed
symmetric operator $\tilde X$ affiliated with $\M$ such that
$\tilde X\Omega=\psi$, to which the $X_n$ converge on the common core
$\M'\Omega$, so that $\tilde XY\Omega=Y\psi$ for every $Y\in\M'$.

If we are to have $U(a)\M U(-a)\subset\M$ for all $a\geq 0$, then we must
have $U(a)R\subset R$ for all $a\geq 0$.  In addition, the generator $H$ of the
unitary group $U(a)$ must give a derivation $\delta$ of $\M$ by
$\delta(X)=i[H,X]$; however, this derivation will be unbounded, hence defined
only on a dense set, and the problem is to give sufficient conditions for
$\delta$ to generate a semigroup of endomorphisms of $\M$.

In the automorphism case, the relevant result of \cite{BH} can be expressed as
follows:

\bigskip
\noindent{\bf Theorem 4:}~~
{\em
Suppose that $U(a)\Omega=\Omega$ and $U(a)R=R$ for all real $a$, and
that the set $D(\delta)=\[X|X\in\M, i{[H,X]}\in\M]$ is such that
$D(\delta)\Omega$ is a core for $H$.  Then $U(a)\M U(-a)=\M$ for all real $a$.
}
\bigskip

The replacement for $U(a)R=R$ in the endomorphism case is clearly
$U(a)R\subset R$ for all $a\geq 0$; we must also find a sufficient replacement
for the core condition.  Let
\begin{equation}
\M_\epsilon=\[X|U(a)XU(-a)\in\M\,\,\hbox{\rm for all}\,\,
0\leq a\leq\epsilon],
\end{equation}
and let $R_\epsilon=\overline{\M_\epsilon^\sa\Omega}$; then
$\M_\epsilon\supset\M_{\epsilon'}$ and $R_\epsilon\supset R_{\epsilon'}$
whenever $\epsilon'\geq\epsilon$.  In addition, let
\begin{equation}
\M_+=\bigcup_{\epsilon>0}\M_\epsilon\quad\hbox{\rm and}\quad
R_+=\bigcup_{\epsilon>0}R_\epsilon.
\end{equation}
Then $\M_\epsilon$ contains those elements $X$ of $\M$ for which the
differential equation $X(t)'=\delta(X(t))$, $X(0)=X$ in the Banach space $\M$
has a solution curve of length at least $\epsilon$; likewise, $\M_+$ contains
those for which there is a solution curve of any positive length.  Conditions
on $\M_\epsilon$ and $\M_+$ can thus be regarded as local existence conditions
for this differential equation; we will use criteria of this sort to control
the behavior of $\delta$.

\bigskip
\noindent{\bf Theorem 5:}~~
{\em
With the above notation and assumptions, suppose that $U(a)\Omega=\Omega$ and
$U(a)R\subset R$ for all $a\geq 0$, and that for some $\epsilon>0$, $\Omega$ is
cyclic for $\M_\epsilon$, i.e., $R_\epsilon+iR_\epsilon$ is dense in $\H$.
Then $U(a)\M U(-a)\subset\M$ for all $a\geq 0$.
}

\smallskip
\noindent{\bf Proof of Theorem 5:}~~
{
It will suffice to show that $U(a)\M'U(-a)\supset\M'$ for all $a\geq 0$;
we have from our assumptions that $U(a)R'\supset R'$ for all $a\geq 0$.  Let us
pick $a$ such that $0\leq a\leq\epsilon$, so that
$\M_\epsilon\subset \M\cap U(-a)\M U(a)$.  Let $Y$ be a self-adjoint element of
$\M'$; then $Y\Omega\in R'\subset U(a)R'$, so that $U(-a)Y\Omega\in R'$.  Then
there is a closed symmetric operator $\tilde Y$ affiliated with $\M'$ such that
$\tilde Y\Omega=U(-a)Y\Omega$, defined on the core $\M\Omega$ by
$\tilde YX\Omega=XU(-a)Y\Omega$ for every $X\in\M$.  $\tilde Y$ therefore
agrees with the bounded operator $U(-a)YU(a)$ on the dense set
$\M_\epsilon\Omega$, from which it follows that $\tilde Y$ is in fact
bounded and equal to $U(-a)YU(a)$, so that $Y\in U(a)\M'U(-a)$ and
$U(a)\M'U(-a)\supset\M'$.  This is so for all $0\leq a\leq\epsilon$, hence by
the semigroup property for all $a\geq 0$.
}
\bigskip

\noindent{\bf Theorem 6:}~~
{\em
With the above notation and assumptions, suppose that
$U(a)\Omega=\Omega$ and $U(a)R\subset R$ for all $a\geq 0$, that
$\Delta^{it}U(a)\Delta^{-it}=U(e^{-2\pi t}a)$ for all real $a$,$t$, and
that $\Omega$ is cyclic for $\M_+$, i.e., $R_++iR_+$ is dense in $\H$. Then
$U(a)\M U(-a)\subset\M$ for all $a\geq 0$.
}
\smallskip

\noindent{\bf Proof of Theorem 6:}~~
{
Notice that $\Delta^{it}R_\epsilon=R_{e^{-2\pi t}\epsilon}$, so for any
$\epsilon>0$, we have $R_+=\cup_{t\geq 0}\lbrace\Delta^{it}R_\epsilon\rbrace$.
By assumption, for any $\psi\in\H$, there is some $\epsilon>0$ and some
$\phi\in R_\epsilon+iR_\epsilon$ such that $\<\psi|\phi>\neq 0$.  Thus given
any $\epsilon>0$, there is some $\phi\in R_\epsilon+iR_\epsilon$ and
some $t\geq 0$ such that $\<\psi|\Delta^{it}\phi>\neq 0$.  But
$\phi\in R+iR=D(\Delta^{1/2})$, so that $\phi$ is an analytic vector for
$\Delta^{it}$ in the strip $-1/2\leq\Im t\leq 0$.  Thus
$\<\psi|\Delta^{it}\phi>$ is the boundary value of a function analytic in $t$
on that strip, and cannot vanish for all $t\leq 0$.  So
$R_\epsilon+iR_\epsilon=\cup_{t\leq 0}
\lbrace\Delta^{it}(R_\epsilon+iR_\epsilon)\rbrace$ is dense in
$\H$ already, and Theorem 5 applies.
}
\goodbreak\bigskip

The condition that $R_++iR_+$ be dense will be referred to as the local
existence condition of Theorem 6; the condition of Theorem 5 is a uniform
version of it.  In specific cases, for example those involving perturbations
of known endomorphism semigroups, we might hope to establish local
existence conditions of these sorts by means of fixed point theorems and
other standard methods for differential equations.

At this point, it seems worthwhile to present the motivating example for this
discussion, in the simple form of a massive scalar free field in 1+1 spacetime
dimensions.  Let $\h$ be a Hilbert space, the one-particle space, and let
$\H=\exp(\h)$ be a symmetric Fock space constructed over it, whose $n$-particle
subspace $\H^n$ is the $n$-fold symmetric tensor product of $\h$ with itself;
$\H^0$ is a one-dimensional Hilbert space, identified with the complex
multiples of the vacuum $\Omega$.  The vectors of $\H$ we will index by the
exponential map for vectors
\begin{equation}
\exp(f)=\sum^\infty_{n=0}{1\over{\sqrt{n!}}}f^{(n)}
\qquad\hbox{ for $f\in\h$,}
\end{equation}
where $f^{(n)}\in\H^n$ is the $n$-fold tensor product of $f$ with itself,
so that $\<\exp(f)|\exp(g)>=\exp(\<f|g>)$; then for any $f\in\h$ we can define
the unitary Weyl operator $\w(f)$ by
\begin{equation}
\w(f)\exp(g)=e^{-{1\over 2}\|f\|^2-\<f| g>}\exp(f+g).
\end{equation}
If $u$ is a unitary operator on $\h$, then its multiplicative second
quantization $U$ given by $U\exp(f)=\exp(uf)$ will be a unitary operator on
$\H$; if $a$ is a self-adjoint operator on $\h$, then its additive second
quantization $A$, the generator of the multiplicative second quantization of
$u(t)=e^{ita}$, will be a self-adjoint operator on $\H$.  The additive second
quantization of the identity is an operator $N$, the number operator, which has
the eigenvalue $n$ on $\H^n$.  Then for any $f\in\h$,
$D_S=\cap_{n=1}^\infty D(N^n)$ will be a core for the generator $\phi(f)$ of
$\w(tf)$, such that $\phi(f)D_S\subset D_S$.

Let $\h$, realized as $\h=\L^2(\RR,d\k)$ with $\nu=-i\partial_\k$, carry the
simplest non-trivial representation of the relations of Theorem 3, realized
by $\delta^{it}=e^{2\pi i\k}$, $u(a)=e^{iae^\nu}$, and
$jf(\k)=\overline{f(-\k)}$.  The multiplicative second quantizations
$\Delta^{it}$, $U(a)$, and $J$ of $\delta^{it}$, $u(a)$, and $j$ respectively
also satisfy the relations of Theorem 3, with a one-dimensional trivial space
(the vacuum space $\H^0$) and a non-trivial representation of infinite
multiplicity.

Then we will let $\M$ (intended to correspond to the right wedge $W_R$) be the
von Neumann algebra generated by $\w(f)$ for all $f$ in the real linear space
\begin{equation}
r=\[f|j\delta^{1/2}f=f]=\[f(\k)=g(\k)+e^{-\pi\k}\overline{g(-\k)}
|g\in D(e^{\pi\k})].
\end{equation}
It can be shown \cite{Ri} that $\M'$ is generated by $\w(f)$ for all $f\in r'$,
where
\begin{equation}
r'=\[f|j\delta^{-1/2}f=f]=\[f(\k)=g(\k)+e^{\pi\k}\overline{g(-\k)}
|g\in D(e^{-\pi\k})],
\end{equation}
and that $J$ and $\Delta^{it}$ give the modular conjugation and automorphisms
for $\M$ with respect to $\Omega$.  Then the unbounded operators $\phi(f)$ for
$f\in r$ will be self-adjoint and affiliated with $\M$, and in fact will
generate $\M$.  It is easy to see that $u(a)r\subset r$ for $a\geq 0$, and thus
that $U(a)\M\subset\M$ for $a\geq 0$.  Let us ignore this for the moment,
however, and proceed to apply Theorems 5 and 6.

Clearly $h_{\lambda,\rho}=\lambda e^\n+\rho e^{-\n}$ is an unbounded
self-adjoint operator on $\h$ for each $(\lambda,\rho)\in\RR^2\setminus(0,0)$;
we may then define the self-adjoint operator $H_{\lambda,\rho}$ as the
additive second quantization of $h_{\lambda,\rho}$, or alternatively by
$i[H_{\lambda,\rho},\phi(f)]=\phi(ih_{\lambda,\rho}f)$.  Then let
\begin{equation}
r_1=\[f(\k)=(\hat g(\sinh\n)+i\cosh\n\,\hat h(\sinh\n))
\check{\phantom{l}}\,|
g,h\,\,\hbox{\rm real and supported in}\,\,[1,\infty)],
\end{equation}
where $\hat{}$ and $\check{}$ represent direct and inverse Fourier transforms.
It can be shown that for every $(\lambda,\rho)\in\RR^2$, there is some
$\epsilon$ such that for every $f\in r_1$, $\phi(f)$ is affiliated with
$\M_\epsilon$ with respect to $H_{\lambda,\rho}$.  Furthermore $r_1+ir_1$ is
dense in $\h$.  It follows that for every $H_{\lambda,\rho}$, the local
existence condition of Theorem 6 and the uniform local existence condition
of Theorem 5 both hold.

\smallskip
Then for $(\lambda,\rho)\in\RR^2\setminus(0,0)$, we have the following:\\
(i) $H_{\lambda,\rho}$ is positive if and only if $\lambda$ and $\rho$ are
both non-negative;\\
(ii) $\Delta^{it}H_{\lambda,\rho}\Delta^{-it}=
H_{e^{-2\pi t}\lambda,e^{2\pi t}\rho}$, and
$JH_{\lambda,\rho}J=H_{-\lambda,-\rho}$;\\
(iii) $H_{\lambda,\rho}$ generates a one-parameter semigroup of endomorphisms
of $\M$ if and only if $\lambda$ and $-\rho$ are both non-negative; and\\
(iv) $H_{\lambda,\rho}$ generates a one-parameter semigroup of endomorphisms of
$\M'$ if and only if $-\lambda$ and $\rho$ are both non-negative.\\
For mass $m$, $H_{m/2,m/2}$ is the Hamiltonian, $H_{m/2,-m/2}$ the momentum;
$U(a)=e^{iaH_{m,0}}$.
\bigskip

Of course, this is a very simple example, in which it is easy to compute
the effects of the $U(a)$.  In more complicated cases, Theorems 5 and 6 could
perhaps be applied to greater effect.  However, their conditions may well be
more restrictive than is necessary; one might conjecture that the local
existence conditions could be replaced by conditions purely on
$D(\delta)$---for example, as in \cite{BH}, by the condition that
$D(\delta)\Omega$ be a core for $H$.

\section{Lightlike Translations}

Let us return to the situation described in the introduction, and consider
again the conditions (a)--(d).  We know already that (a) and (b) each follow
from the remaining three conditions; we have now to consider (c) and (d).
One branch is available immediately:  suppose that (a)
is satisfied, but $U(a)\Omega$ is not known.  Then
\begin{equation}
\ip<\Omega|U(a)\Omega\vphantom{\Delta^{it}}>
=\ip<\Omega|\Delta^{it}U(a)\Delta^{-it}\Omega>=
\ip<\Omega|U(e^{-2\pi t}a)\Omega\vphantom{\Delta^{it}}>
\end{equation}
is independent of $t$, and hence must be a constant for all $a>0$ and for all
$a<0$.  Taking the limit as $t\rightarrow\infty$, these constants must both be
$1$; but since $U(a)\Omega$ is a unit vector, it must therefore equal $\Omega$
for all $a$.  Thus (a) alone implies (d).  With this out of the way, we proceed
to combine the results of Sections II and III:

\bigskip\goodbreak
\noindent{\bf Theorem 7:}~~
{\em
If $H$, the generator of $U(a)$, is positive and annihilates the vacuum, and if
the local existence condition of Theorem 6 holds, then the isotony relation
$U(a)\M U(-a)\subset\M$ for all $a\geq 0$ (and thus $U(a)\M'U(-a)\subset\M'$
for
all $a\leq 0$) holds if and only if the covariance relations hold in the
form
\begin{equation}
\Delta^{it}U(a)\Delta^{-it}=U(e^{-2\pi t}a)\quad\hbox{and}\quad JU(a)J=U(-a)
\end{equation}
for all real $a$ and $t$.  Likewise the covariance and isotony relations
together imply the positivity of $H$.
}

\noindent{\bf Remarks:}~~
{
Although this is in some respects similar to the case of automorphisms, there
are a number of significant differences.  For example, if $H$ were positive in
the automorphism case, then it would be affiliated with $\M$, and since it
annihilates $\Omega$, it would have to vanish; here, however, $H$ can be
positive and non-trivial.
}

\smallskip
\noindent{\bf Proof of Theorem 7:}~~
{
Theorem 6 allows us to reduce this to a question about the relations
between $U(a)$ and $R$:  it will suffice for the first part of the theorem to
show that $U(a)R\subset R$ for all $a\geq 0$ if and only if covariance holds.
But this follows from Theorem 3, since (a${}'$) of Theorem 3 corresponds to
the covariance relations, and (c${}'$) of Theorem 3 is equivalent to the
condition that $U(a)R\subset R$ for all $a\geq 0$.  The second part of the
theorem also follows from Theorem 3, without the aid of Theorem 6.
}
\bigskip

Corresponding results for the backward lightlike translations
$W(a)=T(a(\hat x_1-\hat x_0))$ can be derived by exchanging $\M$ and $\M'$,
and replacing $a$ by $-a$ in the above.  $W(a)$ should have a negative
generator, and should satisfy Lorentz and TCP covariance in the form
$\Delta^{it}W(a)\Delta^{-it}=W(e^{2\pi t}a)$ and $JW(a)J=W(-a)$.  With
these substitutions, the corresponding theorem obtains.  The situation
for the intermediate case, the spacelike translations $T(a\hat x)$ taking $W_R$
into itself, is somewhat more complicated, although not essentially different:
Theorem 5 still holds, but now the relations between the generator (which
now in some frame of reference is the momentum) and the modular operators are
no longer so simple.  We must just show that
$J\Delta^{1/2}T(a\hat x)\supset T(a\hat x)J\Delta^{1/2}$ for all $a\geq 0$.
For example, if $U(a)$ satisfies the conditions of Theorem 7, and $W(a)$ the
corresponding requirements for a backward lightlike translation, and if $U(a)$
and $W(b)$ commute for all $a$ and $b$, then $U(\lambda a)W(-\rho a)$ gives an
endomorphism semigroup of this intermediate type for any $\lambda,-\rho>0$.
This is just the situation described in the example at the end of Section III.
Combining Theorem 7 with these remarks, we have the following omnibus
theorem, as advertised:

\bigskip
\noindent{\bf Theorem 8:}~~
{\em Let $\M$ be a von Neumann algebra acting on a Hilbert space $\H$, which
with its commutant $\M'$ has a separating and cyclic vector $\Omega$.  Given a
strongly continuous one-parameter group $U(a)$ of unitary operators on $\H$,
for which the local existence condition of Theorem 6 holds, any three of the
following conditions imply the fourth:\\
(a) $\Delta^{it}U(a)\Delta^{-it}=U(e^{-2\pi t}a)$ and $JU(a)J=U(-a)$ for
all real $a$ and $t$;\\
(b) the generator $H$ of the $U(a)$ is positive;\\
(c) $U(a)\M U(-a)\subset\M$ for all $a\geq 0$;\\
(d) $U(a)\Omega=\Omega$ for all $a$.\\
Likewise, any three of the following conditions imply the fourth:\\
(a${}'$) $\Delta^{it}U(a)\Delta^{-it}=U(e^{2\pi t}a)$ and $JU(a)J=U(-a)$ for
all real $a$ and $t$;\\
(b${}'$) the generator $H$ of the $U(a)$ is negative;\\
(c) $U(a)\M U(-a)\subset\M$ for all $a\geq 0$;\\
(d) $U(a)\Omega=\Omega$ for all $a$.\\
In addition, the first part of either (a) or (a${}'$) implies (d); if the
first part of (a) holds, then (c) implies (b), or if all of (a) holds, then
(b) and (c) are equivalent.  Likewise if the first part of (a${}'$) holds, then
(c) implies (b${}'$), or if all of (a${}'$) holds, then (b${}'$) and (c) are
equivalent.  Otherwise, no two of these conditions imply any other.
}

\section{Local Algebras}

Let us now extend the example of Section III to a massive scalar free field
in $d+1$ spacetime dimensions where $d>1$.  This will provide us with an
opportunity to prove a result comparable to that of Theorem 8, but without the
assistance of Theorem 6.

Since the $H_{\lambda,\rho}$ are intended to correspond to momenta, and since
the different components of momentum must commute, the proper way to extend
to higher dimensions is as follows:  take $\h=\L^2(\RR,p_0^{-1}d^dp)$,
where $p_0=\sqrt{m^2+p^2}$.  On it represent the $d+1$-dimensional Poincar\'e
group and TCP transformation:  the translations by
$t(x)=e^{i(x_0p_0-x\cdot p)}$; the rotations by rotations of $p$;
the velocity transformations in the $\hat x_i$ direction by
$p_i\rightarrow p_i\cosh\lambda+p_0\sinh\lambda$; and the PCT
transformation by $\theta f(p)=\overline{f(p)}$.  Then these all have
multiplicative second quantizations on the symmetric Fock space $\H$, which
give a representation of the Poincar\'e group for which the vacuum is the only
invariant vector.

Define $\n_i$ by $p_i=\sqrt{p_0^2-p_i^2}\sinh\n_i$ (so that
$p_0=\sqrt{p_0^2-p_i^2}\cosh\n_i$ and $e^{\n_i}=p_0+p_i$) and $\k_i$ by
$\k_i=i\partial_{\n_i}$ (the derivative being taken with $p_{i'}$ fixed for
$i'\neq i$, so that the $\k_i$ so defined do not commute).  Take
$j=\theta\rho$, where $\rho$ is a rotation by the angle $\pi$ about the
$\hat x_1$ axis, and let $\delta=e^{2\pi\k_1}$, where $\k_1$ is to be the
generator of the velocity transformations in the $\hat x_1$ direction.  Then
again let the von Neumann algebra $\A(W_R)$ be that generated by $\w(f)$ for
all $f$ in the real linear space
\begin{equation}
r=\[f|j\delta^{1/2}f=f]=\[f(\k_1,p_i)=g(\k_1,p_i)
+e^{-\pi\k_1}\overline{g(-\k_1,-p_i)}|g\in D(e^{\pi\k_1})],
\end{equation}
where $i$ runs from $2$ to $d$.  From the earlier discussion we see that the
modular conjugation and automorphisms will have the geometric form described in
the introduction.

Any wedge $W$ whose vertex contains the origin is produced from $W_R$ by
some Lorentz transformation; for such wedges we can define corresponding
von Neumann algebras $\A(W)$ by the corresponding transformation of $\A(W_R)$,
and see directly (for example) that $\A(W_L)=\A(W_R)'$.  It is also not
difficult to show that if $W$ is a rotation of $W_R$, then $\A(W)$ is generated
by the two von Neumann algebras $\A(W)\cap\A(W_R)$ and $\A(W)\cap\A(W_L)$, each
of which has the vacuum for a cyclic and separating vector.  In fact, given any
family of such wedges $W_i$ with a nonempty intersection, the $\A(W_i)$ also
have a nonempty intersection, sufficiently large that it has the vacuum as a
cyclic vector (these regions are spacelike cones).  No one of these wedges
contains any other, so there are no isotony relations to establish.  Once we
seek to add the translations, however, we must establish the isotony relation
$T(x)\A(W)T(-x)\subset\A(W)$ for all $x$ such that $x+W\subset W$.  For the
translations given as above, this is immediate; however, as before we wish to
ignore this in the interests of generality.

\bigskip
\noindent{\bf Theorem 9:}~~
{\em
Let $\H$ be a Hilbert space with a representation of the $d+1$-dimensional
Lorentz group, $d>1$, and the TCP operator $\Theta$, for which the vacuum
$\Omega$ is the only invariant vector; let $\A(W)$ be a family of von Neumann
algebras on $\H$ for all $W$ whose vertices contain the origin, covariant under
the Lorentz group, with $\A(W_L)=\A(W_R)'$, and such that the modular
operators with respect to the vacuum are given geometrically from $\Theta$ and
the representation of the Lorentz group, as described above.  Suppose that if
$W(\neq W_R,W_L)$ is a rotation of $W_R$, then $\A(W)$ is generated by
$\A(W)\cap\A(W_R)$ and $\A(W)\cap\A(W_L)$, each of which has the vacuum for a
cyclic and separating vector.

Let $H_\mu$, $\mu=0\ldots d$, be a Lorentz $d+1$-vector of unbounded
strongly commuting self-adjoint operators on $\H$; suppose they are odd under
$\Theta$, annihilate the vacuum, and have their joint spectrum contained in the
closed forward light cone.  Let
\begin{equation}
\Gamma=\[X|X\in\A(W_R),\,i{[H_\mu,X]}\in\A(W_R)\>\>\hbox{for all $\mu$}],
\end{equation}
and let us suppose further that $\Gamma\Omega$ is a common core for all the
$H_\mu$.

Then the $H_\mu$ are the generators of a representation $T(x)$ of the
$d+1$-dimensional translation group, which together with the given
representation of the Lorentz group forms a representation of the Poincar\'e
group under which the vacuum is the only invariant vector.  We may consistently
define $\A(W)$ for any wedge as a Poincar\'e transformation of $\A(W_R)$, and
the $\A(W)$ so defined satisfy isotony.  If the wedges whose vertices contain
the origin possess the property that for any family $W_i$ with a nonempty
intersection, the intersection of the $\A(W_i)$ has the vacuum for a cyclic
vector, then the same may be said for any family of wedges whose vertices
contain a common point.
}

\smallskip
\noindent{\bf Proof of Theorem 9:}~~
{
The assumptions immediately imply that the $H_\mu$ exponentiate to a
representation $T(x)=e^{ix^\mu H_\mu}$ of the $d+1$-dimensional translation
group.  Since they transform covariantly under the given representation
of the Lorentz group, they must combine with it to produce a representation
of the Poincar\'e group for which the vacuum is the unique invariant vector.
Let $R(W)=\overline{\A(W)^\sa\Omega}$; then by Lorentz and TCP covariance
$T(x)R(W)=R(W)$ whenever $x$ is parallel to the vertex of $W$, and by Theorem 3
we have $T(a(\hat x_1+\hat x_0))R(W_R)\subset R(W_R)$ and
$T(a(\hat x_1-\hat x_0))R(W_R)\subset R(W_R)$ for all $a\geq 0$.  Thus
$T(x)R(W_R)\subset R(W_R)$ for every $x\in\overline{W_R}$.

Let us then show that $T(x)\A(W_R)T(-x)=\A(W_R)$ for all $x$ parallel to the
vertex of $W_R$. Let $\hat x$ be some unit vector thus parallel, and
$U(a)=T(a\hat x)$.  Since $\Gamma\Omega$ is a core for the generator of $U(a)$,
and $U(a)R(W_R)=R(W_R)$ for all $a$, by Theorem 4 we know that
$U(a)\A(W_R)U(-a)=\A(W_R)$ for all $a$.  Then covariance implies that
$T(x)\A(W)T(-x)=\A(W)$ for every wedge $W$ and every $x$ parallel to the
vertex of $W$.

Next let us show that $T(x)\A(W)T(-x)\subset\A(W)$ for every $x$ such that
$x+\overline{W}\subset W$.  Without loss of generality we may suppose $W$ to be
a rotation of $W_R$, let us say the wedge in the $\hat x_2$ direction.  By
making a velocity transformation in the $\hat x_2$ direction, and a translation
in an orthogonal direction, it suffices to show this for the translations
$U(a)=T(a\hat x_2)$.  Now, $\M_1=\A(W)\cap\A(W_R)$ and $\M_2=\A(W)\cap\A(W_L)$
together generate $\A(W)$; the vacuum is cyclic and separating for $\M_1$,
$\M_2$, $\M_3=\A(W')\cap\A(W_R)$, and $\M_4=\A(W')\cap\A(W_L)$, where $W'$ is
the wedge in the $-\hat x_2$ direction.  Furthermore by the result of the
preceding paragraph $U(a)\A(W_R)U(-a)=\A(W_R)$ and $U(a)\A(W_L)U(-a)=\A(W_L)$
for all $a$.  Take $X\in\M_1^\sa$; then for any $a\geq 0$ we have
$U(a)XU(-a)\in \A(W_R)$, but also since $U(a)R(W)\subset R(W)$, we have
$U(a)X\Omega\in R(W)$.  Thus there is a closed symmetric operator $\tilde X$
affiliated with $\A(W)$, and such that $\tilde X\Omega=U(a)X\Omega$.  But
$\tilde X$ and $U(a)XU(-a)$ agree on the dense set $\M_4\Omega$, from which it
follows that $\tilde X$ is in fact bounded and equal to $U(a)XU(-a)$, and that
$U(a)XU(-a)\in\A(W)$.  The same argument applies if $X\in\M_2^\sa$, using the
dense set $\M_3\Omega$.  Thus $U(a)\M_1 U(-a)\subset\M_1$ and
$U(a)\M_2 U(-a)\subset\M_2$ for all $a\geq 0$, from which it follows that
$U(a)\A(W)U(-a)\subset\A(W)$ for all $a\geq 0$.

Then weak closure implies that $T(x)\A(W)T(-x)\subset\A(W)$ for all
$x$ such that $x+\overline{W}\subset\overline{W}$.  This is equivalent to
isotony; the remaining condition follows by translation from the
corresponding condition for wedges whose common point is the origin.
}
\bigskip

More work is required to show that the resulting family is in fact local:
for example, that for {\it any} family of wedges $W_i$ with a nonempty
intersection, the intersection of the $\A(W_i)$ is nonempty.  Of course, the
free field is a very simple example, in which it is easy to compute
the effects of the $T(x)$.  In more complicated cases, Theorem 9 could
perhaps be applied to greater effect.

\newpage
\noindent{\Large\bf Acknowledgements}
\bigskip

\noindent{I would like to thank E. H. Wichmann for introducing me to these
subjects, and S. Doplicher and the Universit\`a di Roma for their hospitality
while this work was being completed.}

\end{document}